
\input amstex
\documentstyle{amsppt}
\magnification = 1200
\pageheight{7.5in}
\vcorrection{-.12in}
\NoBlackBoxes
\NoRunningHeads
\nologo
\TagsOnRight
\topmatter
\title{Schwarzian Derivative Criteria for Valence of Analytic and 
Harmonic Mappings}
\endtitle
\author{Martin Chuaqui, Peter Duren, and Brad Osgood}
\endauthor

\address Facultad de Matem\'aticas, P. Universidad Cat\'olica de Chile,
Casilla 306, Santiago 22, Chile
\endaddress
\email mchuaqui\@mat.puc.cl
\endemail

\address Department of Mathematics, University of Michigan, Ann Arbor,
Michigan 48109--1109
\endaddress
\email duren\@umich.edu
\endemail
 
\address Department of Electrical Engineering, Stanford University,
Stanford, California 94305
\endaddress
\email osgood\@ee.stanford.edu
\endemail

\thanks The authors are supported by Fondecyt Grant \# 1030589. 
\endthanks

\abstract  For analytic functions in the unit disk, general bounds 
on the Schwarzian derivative in terms of Nehari functions are shown 
to imply uniform local univalence and in some cases finite and bounded 
valence.  Similar results are obtained for the Weierstrass--Enneper 
lifts of planar harmonic mappings to their associated minimal surfaces.  
Finally, certain classes of harmonic mappings are shown to have finite 
Schwarzian norm.
\endabstract


\endtopmatter

\document

\flushpar
{\bf \S 1.  Introduction.}
\smallpagebreak

     The {\it Schwarzian derivative} of an analytic locally univalent 
function $f$ is defined by 
$$
{\Cal S}f = (f''/f')' - \tfrac12(f''/f')^2\,. 
$$
Its most important property is M\"obius invariance: 
${\Cal S}(T\circ f)={\Cal S}f$ for every M\"obius transformation 
$$
T(z) = \frac{az+b}{cz+d}\,,\qquad ad-bc\neq 0\,.
$$
More generally, if $g$ is any function analytic and locally univalent 
on the range of $f$, then 
$$
{\Cal S}(g\circ f) = (({\Cal S}g)\circ f){f'}^2 + {\Cal S}f\,.
$$
As a special case, ${\Cal S}(g\circ T) = (({\Cal S}g)\circ T){T'}^2$, 
since ${\Cal S}T=0$ for every M\"obius transformation $T$.  For an 
arbitrary analytic function $\psi$, the functions $f$ with Schwarzian 
${\Cal S}f=2\psi$ are those of the form $f=w_1/w_2$, where $w_1$ and $w_2$ 
are linearly independent solutions of the linear differential equation 
$w''+ \psi w=0$.  It follows that ${\Cal S}f={\Cal S}g$ implies  
$f=T\circ g$ for some M\"obius transformation $T$.  In particular, 
M\"obius transformations are the only functions with  ${\Cal S}f=0$. 

     Nehari [15] found that certain estimates on the Schwarzian imply 
global univalence.  Specifically, he showed that if $f$ is analytic 
and locally univalent in the unit disk $\Bbb D$ and its Schwarzian 
satisfies either  $|{\Cal S}f(z)|\leq2(1-|z|^2)^{-2}$ or 
$|{\Cal S}f(z)|\leq{\pi}^2/2$ for all $z\in\Bbb D$, then $f$ is 
univalent in $\Bbb D$. Pokornyi [19] then stated, and Nehari [16] 
proved, that the condition $|{\Cal S}f(z)|\leq4(1-|z|^2)^{-1}$ also 
implies univalence.  In fact, Nehari [16] unified all three criteria by 
proving that $f$ is univalent under the general hypothesis 
$|{\Cal S}f(z)|\leq 2p(|z|)$, where $p(x)$ is a positive continuous even 
function with the properties that $(1-x^2)^2p(x)$ is nonincreasing on 
the interval $[0,1)$ and no nontrivial solution $u$ of the differential 
equation $u''+pu=0$ has more than one zero in $(-1,1)$.  Equivalently, 
the last condition can be replaced by the requirement that some 
solution of the differential equation has no zeros in $(-1,1)$.  We will 
refer to such functions $p(x)$ as {\it Nehari functions}.  It can be 
seen that the functions $p(x)=(1-x^2)^{-2}$, $p(x)={\pi}^2/4$, and 
$p(x)=2(1-x^2)^{-1}$ are all Nehari functions; the nonvanishing solutions 
are respectively $u=\sqrt{1-x^2}$, $u=\cos(\pi x/2)$, and $u=1-x^2$.  

     The constant 2 in Nehari's general criterion is sharp.  In many important  
cases, when $C>2$ the inequality $|{\Cal S}f(z)|\leq C\,p(|z|)$ 
admits nonunivalent functions $f$. For $C>2$ and $p(x)=(1-x^2)^{-2}$, 
Hille [11] gave an explicit example of an analytic function $f$  
satisfying $|{\Cal S}f(z)|\leq C(1-|z|^2)^{-2}$, yet having infinite  
valence in the unit disk.  However, this condition does imply {\it uniform} 
local univalence in the sense that any two points $z_1$ and $z_2$ where 
$f(z_1)=f(z_2)$ are separated by a certain minimum distance in the 
hyperbolic metric.  Further details and generalizations will appear in 
the next section.

     In previous work [3] we introduced a notion of Schwarzian derivative 
for harmonic mappings, or complex-valued harmonic functions.  In [4] 
we generalized Nehari's univalence criterion to the (conformal) lift 
$\widetilde{f}$ of a harmonic mapping $f$ to its associated minimal 
surface $\Sigma$.  Under an inequality of the form 
$$
|{\Cal S}f(z)| + e^{2\sigma(z)}|K(\widetilde{f}(z))| \leq 2p(|z|)\,, 
\qquad z\in\Bbb D\,, \tag 1
$$
for some Nehari function $p$, we showed that $\widetilde{f}$ is univalent 
in $\Bbb D$.  Here $K$ denotes the Gauss curvature of $\Sigma$ and 
$e^\sigma$ is the relevant conformal factor.  When $f$ is analytic, the 
surface $\Sigma$ is the complex plane and $\widetilde{f}=f$, so that 
$K=0$ and the condition (1) reduces to Nehari's criterion.  

     In the present paper we examine the consequences of a weaker 
inequality in which the right-hand side of (1) is replaced by $C\,p(|z|)$ 
for some constant $C>2$.  We begin by discussing the implications for 
analytic functions, where the valence of $f$ depends crucially on the 
growth of $p(x)$ as $x\to1$.  This sets the stage for a more general study  
of the valence of harmonic lifts.  The paper concludes with some remarks 
about Schwarzian norms of harmonic mappings.

\bigpagebreak
\flushpar
{\bf \S 2.  Valence of analytic functions.}
\smallpagebreak

     The {\it valence} of a function $f$ analytic in $\Bbb D$ is defined by 
$\sup_{w\in{\Bbb C}} n(f,w)$, where $n(f,w)\leq\infty$ is the number of 
points $z\in\Bbb D$ for which $f(z)=w$.
     
     For points $\alpha, \beta \in \Bbb D$, the hyperbolic metric is  
$$
d(\alpha, \beta) = \frac12 \, \log\frac{1 + \rho(\alpha, \beta)}
{1 - \rho(\alpha, \beta)}\,, \qquad \text{where} \quad 
\rho(\alpha,\beta)=\left|\frac{\alpha-\beta}{1-\overline{\alpha}\beta}
\right|\,.
$$
Our point of departure is the following known theorem. 
\proclaim{Theorem A}  Let $f$ be analytic and locally univalent, and 
suppose that  
$$
|{\Cal S}f(z)|\leq \frac{2(1+\delta^2)}{(1-|z|^2)^2}\,, 
\qquad z \in \Bbb D\,, 
$$
for some constant $\delta>0$.  Then $d(\alpha,\beta)\geq
\pi/\delta$ for any pair of points $\alpha,\beta\in\Bbb D$ where 
$f(\alpha)=f(\beta)$.  The lower bound is best possible.
\endproclaim

     This result is essentialy due to B. Schwarz [21] and was further developed 
by Beesack and Schwarz [2].  It was rediscovered independently by Minda 
[14] and Overholt [18].  Here is a  relatively simple proof, which will serve as a 
model for arguments given later in this paper.
\demo{Proof} Suppose $f(\alpha)=f(\beta)$ for $\alpha\neq\beta$ and 
consider the function $g=f\circ T$, where $T$ is a M\"obius 
self-mapping of the disk chosen so that $T(0)=\alpha$ and 
$T(b)=\beta$ for some $b>0$.  Then $g(0)=g(b)$ and 
$d(\alpha,\beta)=d(0,b)$ by the M\"obius invariance of the 
hyperbolic metric. Combining the hypothesis on ${\Cal S}f$ with the 
relation ${\Cal S}g={((\Cal S}f)\circ T){T'}^2$, we find after a 
calculation that 
$$
|{\Cal S}g(\zeta)| \leq \frac{2(1+\delta^2)}{(1-|\zeta|^2)^2}\,, 
\qquad \zeta\in\Bbb D\,.
$$
In other words, the bound on the Schwarzian is M\"obius invariant.

     Now let $q=\frac12{\Cal S}g$ and observe that the condition 
$g(0)=g(b)$ implies that some solution of the differential equation 
$u''+qu=0$ has the property  $u(0)=u(b)=0$.  Without loss of 
generality, we may assume that $u(x)\neq0$ for $0<x<b$.  Let $v=|u|$, 
or $v^2=u\overline{u}$, so that 
$$
vv' = \text{Re}\{\overline{u}u'\} \qquad \text{and} \qquad 
|v'(x)| \leq |u'(x)|\,.
$$
Differentiating the first of these relations and applying the second, 
we find
$$
v(x)v''(x) \geq \text{Re}\bigl\{\overline{u(x)}u''(x)\bigr\} 
= - |u(x)|^2\text{Re}\bigl\{q(x)\bigr\}\,,
$$
which implies that $v''(x)+|q(x)|v(x)\geq0$.  Since $v(x)>0$ in the 
interval $(0,b)$, it follows that $v$ satisfies the differential equation 
$v''+hv=0$, where 
$$
h(x) = - \frac{v''(x)}{v(x)} \leq |q(x)| = \tfrac12 |{\Cal S}g(x)| 
\leq \frac{1+\delta^2}{(1-x^2)^2}
$$
for $0<x<b$.  By the Sturm comparison theorem, if $w(x)$ is a solution 
of the differential equation 
$$
w''(x) + \frac{1+\delta^2}{(1-x^2)^2}\,w(x) = 0 \,,  \tag2
$$
with $w(0)=0$, then $w(x)$ also vanishes at a point $c\in(0,b]$.  
But the solutions of (2) with $w(0)=0$ are constant multiples of 
$$
w(x) = \sqrt{1-x^2}\,\sin\left(\frac{\delta}{2}
\log\frac{1+x}{1-x}\right)\,. \tag3
$$
(See Kamke [12], p. 492, eq. 2.369.)  Thus $w(c)=0$ for 
$$
\frac{\delta}{2}\log \frac{1+c}{1-c} = \pi\,, \qquad 
\text{or} \ \  d(0,c) = \frac\pi\delta\,. 
$$
Consequently,  
$$
d(\alpha,\beta) = d(0,b) \geq d(0,c) =  \frac{\pi}{\delta}\,.
$$
As $\delta\to0$ this yields Nehari's univalence criterion with  
$p(x)=(1-x^2)^{-2}$. 

     To see that the bound is sharp, consider the example 
$$
f(z) = \left(\frac{1+z}{1-z}\right)^{i\delta}, \qquad f(0)=1\,,
$$
given by Hille [10].  A calculation shows that  
$$
{\Cal S}f(z) = \frac{2(1+\delta^2)}{(1-z^2)^2}\,,\quad \text{so that} 
\quad |{\Cal S}f(z)| \leq  \frac{2(1+\delta^2)}{(1-|z|^2)^2}\,.
$$
The function $f$ has infinite valence in the disk.  For instance, 
$f(z)=1$ whenever
$$
\delta \log\frac{1+z}{1-z} = 2n\pi\,, \qquad n=0,1,2,\dots\,.
$$
Defining $x_n$ by the equation 
$$
\delta \log\frac{1+x_n}{1-x_n} = 2n\pi\,, 
$$
we see that $0=x_0<x_1<x_2< \cdots < 1$ and $d(x_n,x_{n+1})=\pi/\delta$.  
Thus the lower bound $\pi/\delta$ is best possible. This concludes the 
proof. 
\enddemo 

     We now turn to the more general case where 
$$
|{\Cal S}f(z)| \leq C\,p(|z|)
$$
for some Nehari function $p$ and some constant $C>2$.  The valence of $f$ 
then depends on the size of $C$ and the growth of $p(x)$ as $x\to1$.  
We begin with a general discussion of Nehari functions.  

     Since $(1-x^2)^2p(x)$ is positive and nonincreasing on the interval 
$(0,1)$, the limit 
$$
\mu = \lim_{x\to1-} (1-x^2)^2p(x)  \tag4
$$
exists and $\mu\geq 0$.  Observe first that $\mu\leq1$.  Indeed, if 
$\mu>1$ then 
$$
p(x) > \frac{\frac12(1+\mu)}{(1-x^2)^2}
$$
in some interval $x_0<x<1$, and the Sturm comparison theorem shows 
that the solutions of the differential equation $u''+pu=0$ have 
infinitely many zeros in $(-1,1)$, which is not possible because 
$p$ is a Nehari function.

     If $\mu=1$, we claim that $p(x)$ must be the function $(1-x^2)^{-2}$.  
The preceding argument does not quite capture this result.  The proof 
will appeal instead to the following general extension of the Sturm theory, 
implicit in previous papers by Chuaqui and Osgood [7], Chuaqui and 
Gevirtz [6] and Chuaqui, Duren, and Osgood [4].  
\proclaim{Relative Convexity Lemma}  Let $P$ and $Q$ be continuous 
functions on a finite interval $[a,b]$, with $Q(x)\leq P(x)$.  Let $u$ 
and $v$ be solutions to the respective differential equations 
$u''+Pu=0$ and $v''+Qv=0$, with $u(a)=v(a)=1$ and $u'(a)=v'(a)=0$.  
Suppose that $u(x)>0$ and $v(x)>0$ in $[a,b)$ and define the function 
$$
F(x) = \int_a^x v(t)^{-2}\, dt\,, \qquad a\leq x < b\,.
$$
Then $F$ is continuous and increasing on $[a,b)$, and it maps this 
interval onto an interval $[0,L)$, where $0<L\leq\infty$.  Let $G$ denote 
the inverse of $F$.  Then the function 
$$
w(y) = \frac{u(G(y))}{v(G(y))} 
$$
is concave on the interval $[0,L)$.  
\endproclaim
\demo{Proof} We are to show that $w''(y)\leq 0$.  Differentiation gives 
$$
w' = \frac{vu' - uv'}{v^2}\,G' = vu' - uv'\,, 
$$
since $G'(y)=1/F'(G(y))=v(G(y))^2$.  Another differentiation gives 
$$
w'' = (vu'' - uv'')G' = uv(Q-P)v^2 = (Q-P)v^4w\,.
$$
Since $Q(x)\leq P(x)$ and $w(y)>0$ by hypothesis, this shows that 
$w''(y)\leq0$ for $0\leq y<L$.
\enddemo
 
     With the lemma in hand, we now show that $p(x)=(1-x^2)^{-2}$ if 
$\mu=1$.  Because $(1-x^2)^2p(x)$ decreases to 1, it is clear that 
$p(x)\geq(1-x^2)^{-2}$.  Let $q(x)=(1-x^2)^{-2}$ and note that 
$v(x)=\sqrt{1-x^2}$ is the solution of $v''+qv=0$ on $[0,1)$ with 
$v(0)=1$ and $v'(0)=0$.  Thus
$$
F(x) = \int_0^x v(t)^{-2}\,dt = \frac12 \log \frac{1+x}{1-x} 
\to \ \infty \qquad \text{as} \ \ x\to1\,, 
$$
and the inverse function $G(y)$ is defined on the interval 
$0\leq y<\infty$.  If $u$ is the solution of $u''+pu=0$ with $u(0)=1$ 
and $u'(0)=0$, the relative convexity lemma says that 
$w(y)=u(G(y))/v(G(y))$ is concave on $[0,\infty)$, with $w(0)=1$ and 
$w'(0)=0$.  If $p(x)\not\equiv q(x)$, then $w(y)\not\equiv1$ and so 
$0<w(y_0)<1$ and $w'(y_0)<0$ for some point $y_0\in(0,\infty)$.  Since 
the curve $w(y)$ is concave, it lies below its tangent line constructed 
at $y_0$, and this line must cross the $y$-axis.  Therefore, $w(y_1)=0$ 
for some point $y_1\in(y_0,\infty)$, which implies that $u(x_1)=0$ for 
some point $x_1\in(0,1)$.  But this is not possible because by definition 
of a Nehari function the solution $u$ cannot vanish more than once in 
$(-1,1)$.  However, $u$ is an even function since $p$ is even, so 
$u(x_1)=0$ would imply $u(-x_1)=0$.  This contradiction shows that 
$p(x)\equiv q(x)=(1-x^2)^{-2}$ in $[0,1)$, hence in $(-1,1)$.

     Consequently, we may suppose that $0\leq\mu<1$.  It is relevant 
to note that each such value of $\mu$ actually arises from some 
Nehari function.  For a parameter $t$ in the interval $1<t<2$,
$$
p(x) = \frac{t(1-(t-1)x^2)}{(1-x^2)^2}\,,  \tag5
$$
is a Nehari function with $\mu=t(2-t)$.  Indeed, it can be checked that 
the function $u(x)=(1-x^2)^{t/2}$ solves the differential equation 
$u''+pu=0$ and has no zeros on $(-1,1)$.  The other required properties 
of a Nehari function are easily verified, so this shows that each 
$\mu\in(0,1)$ arises from a Nehari function.  The Nehari functions 
$p(x)=\pi^2/4$ and $p(x)=2(1-x^2)^{-1}$, mentioned in the introduction, 
have $\mu=0$.

     Incidentally, it can be shown that for $1<t<2$ the function
$$
f(z) = \int_0^z \frac{d\zeta}{(1-\zeta^2)^t} 
$$
has Schwarzian derivative satisfying $|{\Cal S}f(z)| \leq 2 p(|z|)$ for 
all $z\in\Bbb D$\,, where $p$ is the Nehari function defined by (5).  Thus 
by Nehari's theorem, $f$ is univalent in $\Bbb D$\,, a result that Nehari 
[17] obtained by a different method.  

     The following theorem is an analogue of Theorem A for other Nehari 
functions.  
\proclaim{Theorem 1} Let $f$ be analytic and locally univalent in 
$\Bbb D$\,, and suppose that 
$$
|{\Cal S}f(z)| \leq C\,p(|z|)\,, \qquad z\in\Bbb D\,,
$$
for some constant $C>2$ and some Nehari function $p$ with $0\leq\mu<1$, 
where $\mu$ is defined by $(4)$.  Then $f$ is uniformly locally 
univalent with respect to the hyperbolic metric.  If $C\mu<2$, the 
valence of $f$ is finite and has a bound independent of $f$.  
\endproclaim
\proclaim{Corollary} Let $p$ be a Nehari function with $\mu=0$, and 
suppose that 
$$
|{\Cal S}f(z)| \leq C\,p(|z|)\,, \qquad z\in\Bbb D\,,
$$ 
for some constant $C$. Then $f$ has finite valence in $\Bbb D$\,.  
\endproclaim 

     In the statement of the theorem the case $\mu=1$ is excluded 
because the uniform local univalence is then guaranteed by Theorem A.  
The case $C\leq2$ is excluded because the condition 
$|{\Cal S}f(z)| \leq 2\,p(|z|)$ is known to imply that $f$ is univalent 
(Nehari's theorem).  
\demo{Proof of Theorem 1}  Because $\mu<1$ it is clear that 
$p(x)\leq(1-x^2)^{-2}$ for all $x$ in some interval $r<x<1$, and so 
$$
|{\Cal S}f(z)| \leq C\,p(|z|)\leq\frac{C}{(1-|z|^2)^2}\,,\qquad r<|z|<1\,. 
$$
Consequently, $|{\Cal S}f(z)| \leq C_1(1-|z|^2)^{-2}$ for some constant 
$C_1\geq C$ and all $z\in\Bbb D$.  
Setting $C_1=2(1+\delta^2)$ with $\delta>0$, one can infer from 
Theorem A that $d(\alpha,\beta)\geq\pi/\delta$ for any pair of points 
where $f(\alpha)=f(\beta)$.  Thus  $f$ is uniformly locally univalent 
with respect to the hyperbolic metric.  

     If $C\mu<2$, we can say more.  Then in some annulus $r<|z|<1$ 
we have 
$$
|{\Cal S}f(z)| \leq C\,p(|z|) \leq \frac{C_1}{(1-|z|^2)^2} 
$$
for some constant $C_1<2$.  By a theorem of Gehring and Pommerenke 
([10], Theorem 4), this implies that $f$ has finite valence in the annulus 
$r<|z|<1$, and there is a bound on its valence that depends only on $C$ and $p$, 
not on $f$.  In the disk $|z|\leq r$ the uniform local univalence again shows 
that $f$ has finite valence with a bound independent of $f$.  Thus the same 
is true in $\Bbb D$\,, as the theorem asserts.  It should be remarked that a 
result slightly weaker than that of Gehring and Pommerenke is implicit in 
earlier work of B. Schwarz ([21], Theorem 1).  The  result of Schwarz shows 
that the condition $|{\Cal S}f(z)| \leq C(1-|z|^2)^{-2}$ implies $n(f,w)<\infty$ for 
each $w\in\Bbb C$ but does not provide a bound that is uniform in $w$.   
\enddemo

     As applications of Theorem 1, if $|{\Cal S}f(z)| \leq C$ or if 
$|{\Cal S}f(z)| \leq C(1-|z|^2)^{-1}$ in $\Bbb D$ for any constant $C$, 
then $f$ has finite valence.  In either case it would be interesting 
to find a sharp bound on the valence of $f$ in terms of $C$.  If  
$|{\Cal S}f(z)|\leq C$, we can adapt the proof of Theorem A to show 
that any pair of points $\alpha,\beta\in\Bbb D$ where 
$f(\alpha)=f(\beta)$ must satisfy $|\alpha-\beta|\geq\sqrt{2/C}\,\pi$.
(In particular, if $|{\Cal S}f(z)|\leq\pi^2/2$, then    
$|\alpha-\beta|\geq2$ and so $f$ is univalent in $\Bbb D$, in 
accordance with Nehari's theorem.)  The lower bound is established 
by comparing the associated differential equation $u''+({\Cal S}f/2)u=0$ 
with $v''+(C/2)v=0$, whose solutions have zeros with separation 
$\sqrt{2/C}\,\pi$.  An explicit upper bound on the valence of $f$ then 
follows by consideration of the optimal packing of disks of fixed radius 
$\sqrt{1/2C}\,\pi$ in the unit disk.  Details are pursued in a subsequent 
paper [5], which gives some quantitative estimates on the valence.

\bigpagebreak
\flushpar
{\bf \S 3.  Valence of harmonic lifts.}
\smallpagebreak

     We now turn to harmonic mappings and show how the preceding 
theorems on valence of analytic functions can be generalized to 
the lifts of harmonic mappings to their associated minimal surfaces. 
A planar harmonic mapping is a complex-valued harmonic function
$$
f(z) = u(z) + i v(z)\,,\qquad z=x+iy\,,
$$
defined on some domain $\Omega\subset\Bbb C$.  If $\Omega$ is simply 
connected, the mapping has a canonical decomposition $f=h+\overline{g}$, 
where $h$ and $g$ are analytic in $\Omega$ and $g(z_0)=0$ for some specified 
point $z_0\in \Omega$.  The mapping $f$ is locally univalent if and only if 
its Jacobian $|h'|^2 - |g'|^2$ does not vanish.  It is said to be
orientation-preserving if $|h'(z)|>|g'(z)|$ in $\Omega$, or equivalently if
$h'(z)\neq0$ and the dilatation $\omega=g'/h'$ has the property
$|\omega(z)|<1$ in $\Omega$.

     According to the Weierstrass--Enneper formulas, a harmonic mapping 
$f=h+\overline{g}$ with $|h'(z)|+|g'(z)|\neq0$ lifts locally to 
a minimal surface described by conformal parameters if and only if 
its dilatation has the form $\omega=q^2$ for some meromorphic function $q$.  
The Cartesian coordinates $(U,V,W)$ of the surface are then given by
$$
U(z)=\text{Re}\{f(z)\}\,,\quad 
V(z)=\text{Im}\{f(z)\}\,,\quad
W(z)= 2\,\text{Im}\left\{\int_{z_0}^z 
\sqrt{h'(\zeta)g'(\zeta)}\,d\zeta\right\}\,.
$$
We use the notation 
$$
\widetilde{f}(z) = \bigl(U(z),V(z),W(z)\bigr)
$$
for the lifted mapping from $\Omega$ to the minimal surface.  The first 
fundamental form of the surface is $ds^2=\lambda^2|dz|^2$, where the 
conformal metric is 
$$
\lambda=e^\sigma = |h'|+|g'|\,.
$$
The Gauss curvature of the surface at a point $\widetilde{f}(z)$ for 
which $h'(z)\neq0$ is 
$$
K = - \frac{1}{\lambda^2} \Delta (\log \lambda) 
= - \frac{4|q'|^2}{|h'|^2(1+|q|^2)^4}\,,  \tag6
$$  
where $\Delta$ is the Laplacian operator.  Further information about 
harmonic mappings and their relation to minimal surfaces can be found 
in the book [9].
 
   For a harmonic mapping $f=h+\overline{g}$ with $|h'(z)|+|g'(z)|\neq0$, 
whose dilatation is the square of a meromorphic function, we have 
defined [3] the {\it Schwarzian derivative} by the formula 
$$
{\Cal S}f = 2\bigl(\sigma_{zz} - \sigma_z^2\bigr)\,, \tag7
$$
where $\sigma = \log\bigl(|h'|+|g'|\bigr)$ and  
$$
\sigma_z = \frac{\partial\sigma}{\partial z}  
= \frac12 \left(\frac{\partial\sigma}{\partial x} 
- i \frac{\partial\sigma}{\partial y}\right)\,, \qquad z = x+iy\,.
$$
If $f$ is analytic, ${\Cal S}f$ is the classical Schwarzian.  If $f$ 
is harmonic and $\phi$ is analytic, then $f\circ\phi$ is harmonic and 
$$
{\Cal S}(f\circ\phi) = (({\Cal S}f)\circ\phi){\phi'}^2 + {\Cal S}\phi\,,
$$
a generalization of the classical formula.  With $h'(z)\neq0$ and 
$g'/h'=q^2$, a calculation ({\it cf}. [3]) produces the expression
$$
{\Cal S}f = {\Cal S}h +\frac{2\overline{q}}{1+|q|^2}
\left(q'' - q'\,\frac{h''}{h'}\right) -4\left(\frac{q'\overline{q}}{1+|q|^2}
\right)^2\,.\tag8
$$ 

     It must be emphasized that we are not requiring our harmonic 
mappings to be locally univalent.  In other words, the Jacobian need 
not be of constant sign in the domain $\Omega$.  The orientation of 
the mapping may reverse, corresponding to a folding in the associated 
minimal surface.  It is also possible for the minimal surface to 
exhibit several sheets above a point in the $(u,v)$--plane.  Thus the 
lifted mapping $\widetilde{f}$ may be locally or globally univalent 
even when the underlying mapping $f$ is not.   

     In earlier work [4] we obtained the following criterion for the 
lift of a harmonic mapping to be univalent.

\proclaim{Theorem B}  Let $f=h+\overline{g}$ be a harmonic mapping of the 
unit disk, with $e^{\sigma(z)}=|h'(z)|+|g'(z)|\neq0$ and dilatation 
$g'/h'=q^2$ for some meromorphic function $q$.  Let $\widetilde{f}$ 
denote the Weierstrass--Enneper lift of $f$ to a minimal surface with 
Gauss curvature $K=K(\widetilde{f}(z))$ at the point $\widetilde{f}(z)$.  
Suppose that the inequality $(1)$ holds   
for some Nehari function $p$. Then $\widetilde{f}$ is univalent in $\Bbb D$.
\endproclaim

     Examples in [4] show that the univalence may fail if the right-hand 
side of (1) is replaced by $Cp(z)$ for some $C>2$, but we now give 
harmonic analogues of Theorem A and Theorem 1.  Recall that $\mu$ is 
defined by the limit (4).
\proclaim{Theorem 2}  Let $f=h+\overline{g}$ be a harmonic mapping of the 
unit disk, with $e^{\sigma(z)}=|h'(z)|+|g'(z)|\neq0$ and dilatation 
$g'/h'=q^2$ for some meromorphic function $q$.  Let $\widetilde{f}$ denote 
the Weierstrass--Enneper lift of $f$ to a minimal surface with 
Gauss curvature $K=K(\widetilde{f}(z))$ at the point $\widetilde{f}(z)$.  
Suppose that   
$$
|{\Cal S}f(z)| + e^{2\sigma(z)} |K(\widetilde{f}(z))| \leq C\,p(|z|)\,, 
\qquad z\in\Bbb D\,, \tag9
$$
for some Nehari function $p$ and some constant $C=2(1+\delta^2)$, where 
$\delta>0$.  Then $\widetilde{f}$ is uniformly locally univalent. 
If $\mu=1$, so that $p(x)=(1-x^2)^{-2}$, then 
$d(\alpha,\beta)\geq\pi/\delta$ for any pair of distinct points in 
$\Bbb D$ where $\widetilde{f}(\alpha)=\widetilde{f}(\beta)$.  If $C\mu<2$, 
then $\widetilde{f}$ has finite valence that is bounded independently of $f$.
\endproclaim

     The proof will rely on devices similar to those used in the proof 
of Theorem B.  The first is a notion of Schwarzian derivative for curves 
in ${\Bbb R}^n$, due to Ahlfors [1].  For a sufficiently smooth 
mapping $\varphi: (a,b)\to{\Bbb R}^n$, Ahlfors defined 
$$
S_1\varphi = \frac{\langle \varphi',\varphi'''\rangle}{|\varphi'|^2}
- 3\frac{\langle \varphi',\varphi''\rangle^2}{|\varphi'|^4} 
+ \frac32\frac{|\varphi''|^2}{|\varphi'|^2}
$$
and proved its invariance under postcomposition with M\"obius 
transformations.  Chuaqui and Gevirtz [6] recently used Ahlfors' 
Schwarzian to give a univalence criterion for curves.  Here is their 
result:
\proclaim{Theorem C} Let $P(x)$ be continuous on an interval $(a,b)$, 
with the property that no nontrivial solution of the differential 
equation $u''+Pu=0$ has more than one zero.  Let 
$\varphi: (a,b)\to{\Bbb R}^n$ be a curve of class $C^3$ with 
$\varphi'(x)\neq0$ on $(a,b)$.  If $S_1\varphi(x)\leq2P(x)$ on $(a,b)$, 
then $\varphi$ is univalent.  
\endproclaim
\demo{Proof of Theorem 2}  Suppose first that $f$ satisfies (9) with 
$C=2(1+\delta^2)$ and $p(x)=(1-x^2)^{-2}$.  It may be observed that for 
this choice of Nehari function the hypothesis is M\"obius invariant.  More 
precisely, if $T$ is any M\"obius self-mapping of the disk, then 
$F(\zeta)=f(T(\zeta))$ satisfies 
$$
|{\Cal S}F(\zeta)| + e^{2\tau(\zeta)} |K(\widetilde{F}(\zeta))| 
\leq \frac{C}{(1-|\zeta|^2)^2}\,, \qquad \zeta\in\Bbb D\,, \tag10
$$
where $\widetilde{F}=\widetilde{f}(T(\zeta))$ is the lift of the $F$ 
to the same minimal surface, now endowed with the conformal metric 
$e^{\tau(\zeta)}=e^{\sigma(T(\zeta))}|T'(\zeta)|$.  But 
$$
{\Cal S}F(\zeta) = ({\Cal S}f)(T(\zeta))T'(\zeta)^2\,,
$$
by the transformation property of the Schwarzian, so (10) follows from (9) 
via the relation 
$$
\frac{|T'(\zeta)|}{1-|T(\zeta)|^2} = \frac{1}{1-|\zeta|^2}\,.
$$

     Now suppose that $\widetilde{f}(z_1)=\widetilde{f}(z_2)$ for some 
pair of distinct points $z_1, z_2\in\Bbb D$\,.  Let $T$ be the M\"obius 
transformation of $\Bbb D$ onto itself such that $T(0)=z_1$ and $T(\xi)=z_2$ 
for some point $\xi$ in the real interval $(0,1)$.  Let $F=f\circ T$ 
and consider the curve $\varphi(x)=\widetilde{F}(x)$, for which 
$\varphi(0)=\varphi(\xi)$.  A calculation ({\it cf.} [4], Lemma 1) 
shows that
$$
S_1\widetilde{F}(x) \leq \text{Re}\{{\Cal S}F(x)\} 
+ e^{2\tau(x)}|K(\widetilde{F}(x))|\,, \qquad -1<x<1\,.
$$
Consequently, it follows from (10) that 
$$
S_1\varphi(x) \leq \frac{C}{(1-|\zeta|^2)^2}\,, \qquad -1<x<1\,.
$$
This allows an application of Theorem C with 
$P(x)=(1+\delta^2)(1-x^2)^{-2}$.  Recall that the differential equation 
(2) has a solution (3) that does not vanish in an interval $(0,b)$ 
whose endpoints are at hyperbolic distance $d(0,b)=\pi/\delta$.  
In view of the Sturm separation theorem, this implies that no nontrivial 
solution of (2) can vanish more than once in the interval $[0,b)$.  Thus 
Theorem C applies to show that $\varphi$ is univalent in $[0,b)$.  
But $\varphi(0)=\varphi(\xi)$ and so $d(0,\xi)\geq d(0,b)=\pi/\delta$. 
By the M\"obius invariance of the hyperbolic metric, it follows that 
$d(z_1,z_2)\geq\pi/\delta$.  (Strictly speaking, $0\notin(0,b)$ and 
so Theorem C must be applied to a sequence of intervals $(a_n,b_n)$ 
with $a_n<0<b_n$, $d(a_n,b_n)=\pi/\delta$, and $a_n\to0$.)

     For any other Nehari function $p(x)$, the property $\mu\leq1$ implies that 
$C\,p(|z|)\leq C_1(1-|z|^2)^{-2}$ in $\Bbb D$, and so if (9) holds it follows from 
what we have just shown that $\widetilde{f}$ is uniformly locally univalent.

     Suppose next that $C\mu<2$.  To conclude that $\widetilde{f}$ 
has bounded valence, we adapt the argument given by Gehring and Pommerenke 
[10] for analytic functions.  First note that $C\mu<2$ implies $\mu<1$, 
and so by definition of $\mu$ it follows from (9) that 
$$
|{\Cal S}f(z)| + e^{2\sigma(z)} |K(\widetilde{f}(z))| \leq 
\frac{C_1}{(1-|z|^2)^2} \tag11
$$
in some annulus $r<|z|<1$, for some constant $C_1<2$, where the radius 
$r$ depends on $C$ and $p$ but not on $f$.    Now let $\Omega$ 
be an ``almost circular'' domain in the annulus $r<|z|<1$ whose boundary 
consists of a small arc of the unit circle together with an arc of a 
circle tangent to the circle $|z|=r$.  Let $\psi$ be a conformal mapping 
of $\Bbb D$ onto $\Omega$, and define $F=f\circ\psi$.  Then the 
Schwarzian of $F$ is 
$$
{\Cal S}F(\zeta) = ({\Cal S}f)(\psi(\zeta))\psi'(\zeta)^2 + 
{\Cal S}\psi(\zeta)\,,\qquad \zeta\in\Bbb D\,.
$$
Explicit calculation ({\it cf.} [10], p. 238) shows that the Schwarzian norm
$\|{\Cal S}\psi\|$, as defined in (14) below,  can be made arbitrarily small by 
choosing the domain $\Omega$ to be sufficiently circular.  (Intuitively, this makes 
$\psi$ ``almost'' a M\"obius transformation, which has zero Schwarzian.)  
Furthermore, since $\psi$ maps $\Bbb D$ into itself, it satisfies the 
inequality
$$
\frac{|\psi'(\zeta)|}{1-|\psi(\zeta)|^2} \leq \frac{1}{1-|\zeta|^2}\,,
\qquad \zeta\in\Bbb D\,.
$$
But the composite mapping $F=f\circ\psi$ has a lift $\widetilde{F}(\zeta)=
\widetilde{f}(\psi(\zeta))$ to the same minimal surface with conformal 
metric
$$
e^{\tau(\zeta)} = e^{\sigma(\psi(\zeta))}|\psi'(\zeta)|\,.
$$
Combining these relations with the inequality (11), we find that
$$
\align
|{\Cal S}F(\zeta)| + e^{2\tau(\zeta)} |K(\widetilde{F}(\zeta))| &\leq 
\frac{C_1}{(1-|\zeta|^2)^2} + |{\Cal S}\psi(\zeta)| \\
&\leq \frac{C_2}{(1-|\zeta|^2)^2}\,, \qquad \zeta\in\Bbb D\,,
\endalign
$$
where $C_1<C_2<2$, if $\Omega$ is sufficiently circular.  
Invoking Theorem B with $p(x)=(1-x^2)^{-2}$, we conclude that 
$\widetilde{F}$ is univalent in $\Bbb D$.  Therefore, $\widetilde{f}$ 
is univalent in $\Omega$.  

     Now let $\rho=\frac12(1+r)$ and observe that the annulus 
$\rho<|z|<1$ is contained in the union of a finite number of rotated 
copies of $\Omega$.  Since $\widetilde{f}$ is univalent in each of 
these copies of $\Omega$, it follows that $\widetilde{f}$ has finite 
valence in the annulus $\rho<|z|<1$, and the valence has a bound 
independent of $f$.  On the other hand, the uniform local univalence  
shows that $\widetilde{f}$ has finite and bounded valence in each 
closed subdisk $|z|\leq R<1$.  Therefore, $\widetilde{f}$ has finite 
and bounded valence in $\Bbb D$, and the proof of Theorem 2 is complete.
\enddemo

\bigpagebreak
\flushpar
{\bf \S 4.  An example.}
\smallpagebreak

     Our proof of Theorem A essentially incorporates the example 
given by Hille [11] to show that Nehari's original univalence 
criterion is best possible.  We now adapt Hille's example to see  
that the inequality $d(\alpha,\beta)\geq\pi/\delta$ of Theorem 2 is 
best possible for harmonic mappings that lift to a catenoid.  

     The harmonic function $w=f(z)=z+1/\overline{z}$ maps the 
punctured plane $0<|z|<\infty$ onto the doubly covered annular 
region $2\leq|w|<\infty$, and it lifts to the mapping $\widetilde{f}=
(U,V,W)$ of the punctured plane onto the catenoid defined by 
$$
U=(r+1/r)\cos\theta\,,\quad V=(r+1/r)\sin\theta\,,\quad W=2\log r\,,
$$
where $z=re^{i\theta}$.  Consider now the composition $F=f\circ\phi$, 
with
$$
z=\phi(\zeta)= c\left(\frac{1+\zeta}{1-\zeta}\right)^{i\delta} 
= c\,\exp\left\{i\delta\log\frac{1+\zeta}{1-\zeta}\right\}\,, \qquad 
\zeta\in\Bbb D\,.
$$
Here $\delta>0$ and $c>0$ are fixed parameters.  Although $\widetilde{f}$ 
is univalent, the lifted mapping $\widetilde{F}=\widetilde{f}\circ\phi$ 
has infinite valence, since for instance $\phi$ maps the real segment 
$-1<\zeta<1$ onto the infinitely covered circle $|z|=c$.  As in the 
proof of Theorem A, we find that $\widetilde{f}(x_n)=\widetilde{f}(x_{n+1})$ 
for a sequence of points $x_n$ on the real axis at hyperbolic distance 
$d(x_n,x_{n+1})=\pi/\delta$.

     On the other hand, calculations lead to the Schwarzian expressions 
$$
{\Cal S}f(z) = \frac{4\,|z|^2}{z^2(1+|z|^2)^2}\,, \qquad 
{\Cal S}\phi(\zeta) = \frac{2(1+\delta^2)}{(1-\zeta^2)^2}\,,
$$
so the chain rule gives 
$$
{\Cal S}F(\zeta) = \frac{2}{(1-\zeta^2)^2}\left\{1 + \delta^2 
- \frac{8\,\delta^2\,|\phi(\zeta)|^2}{(1+|\phi(\zeta)|^2)^2}\right\}\,.\tag12
$$
By the formula (6) with $h'(z)=1$ and $q(z)=i/z$, the Gauss curvature of 
the catenoid is found to be 
$$
K(\widetilde{f}(z)) = - \frac{4\,|z|^4}{(1+|z|^2)^4}\,,
$$
and so 
$$
e^{2\sigma(z)}K(\widetilde{f}(z)) = - \frac{4}{(1+|z|^2)^2}\,.
$$
Thus 
$$
\aligned
e^{2\tau(\zeta)}|K(\widetilde{F}(\zeta))| &=e^{2\sigma(\phi(\zeta))}
|\phi'(\zeta)|^2 |K(\widetilde{f}(\phi(\zeta)))| \\
&= \frac{4\,|\phi'(\zeta)|^2}{(1+|\phi(\zeta)|^2)^2} 
= \frac{16\,\delta^2\,|\phi(\zeta)|^2}{(1+|\phi(\zeta)|^2)^2\,
|1-\zeta^2|^2}\,. 
\endaligned
\tag13
$$
Finally, the expressions (12) and (13) combine to show that 
$$
|{\Cal S}F(\zeta)| + e^{2\tau(\zeta)}|K(\widetilde{F}(\zeta))| 
= \frac{2(1+\delta^2)}{|1-\zeta^2|^2} \leq 
\frac{2(1+\delta^2)}{(1-|\zeta|^2)^2}
$$
if the constant $c$ is chosen sufficiently small.  In other words,
the harmonic mapping $F$ satisfies the hypothesis of Theorem 2 with 
$C=2(1+\delta^2)$ and Nehari function $p(x)=(1-x^2)^{-2}$, yet 
$d(\alpha,\beta)=\pi/\delta$ at certain points where 
$\widetilde{F}(\alpha)=\widetilde{F}(\beta)$, so the general estimate 
$d(\alpha,\beta)\geq\pi/\delta$ cannot be improved in this case.

\bigpagebreak
\flushpar
{\bf \S 5.  Schwarzian norms of harmonic mappings.}
\smallpagebreak

     Let $f=h+\overline{g}$ be a locally univalent orientation-preserving 
harmonic mapping with dilatation $g'/h'=q^2$, where $q$ is analytic and 
$|q(z)|<1$ in the unit disk $\Bbb D$.  Then the Schwarzian ${\Cal S}f$ 
is defined by (7) and has the expression (8).  Let $\|{\Cal S}f\|$ denote 
the hyperbolic norm of ${\Cal S}f$; that is,
$$
\|{\Cal S}f\| = \sup_{z\in\Bbb D} (1-|z|^2)^2 |{\Cal S}f(z)|\,. \tag14
$$
It should be observed that the Schwarzian norm is M\"obius invariant.  
In other words, $\|{\Cal S}(f\circ T)\| = \|{\Cal S}f\|$ if $T$ is any  
M\"obius self-mapping of the disk. 

     If $f$ is {\it analytic} and univalent in $\Bbb D$, it has long been 
known that $\|{\Cal S}f\|\leq6$, a result due to Kraus [13].  
The bound is sharp, since the Koebe function $k(z)=z/(1-z)^2$ has 
Schwarzian 
$$
{\Cal S}k(z) = - \frac{6}{(1-z^2)^2}\,.
$$
In view of this result for analytic functions, it is natural to ask 
whether every univalent harmonic mapping $f$ has finite Schwarzian norm.  
We have not been able to settle this question, but we can show, with no 
assumption of global univalence, that a harmonic mapping has finite 
Schwarzian norm if and only if its analytic part does.  
\proclaim{Theorem 3}  Let $f=h+\overline{g}$ be a locally univalent 
orientation-preserving harmonic mapping whose dilatation is the square 
of an analytic function in the unit disk.  Then $\|{\Cal S}f\|<\infty$ if 
and only if $\|{\Cal S}h\|<\infty$.  In particular, $\|{\Cal S}f\|<\infty$ 
if $h$ is univalent.
\endproclaim
\demo{Proof} Suppose first that $\|{\Cal S}h\|<\infty$.  By hypothesis, 
$f$ has dilatation $g'/h'=q^2$ for some analytic function $q$ with 
$|q(z)|<1$.  Hence 
$$
\frac{2\,|q(z)|}{1+|q(z)|^2} \leq 1\,, \qquad z\in\Bbb D\,,
$$
and so it follows from (8) that 
$$
|{\Cal S}f(z)| \leq |{\Cal S}h(z)| + |q''(z)| 
+ |q'(z)|\,\left|\frac{h''(z)}{h'(z)}\right| + |q'(z)|^2\,.
$$
By the Schwarz--Pick lemma, 
$$
|q'(z)| \leq \frac{|q'(z)|}{1-|q(z)|^2} \leq \frac{1}{1-|z|^2}\,.
$$
If a function $\phi$ is analytic in the unit disk and 
$|\phi(z)|\leq1/(1-|z|^2)$, then it follows from Cauchy's integral 
formula that $|\phi'(z)|\leq4/(1-|z|^2)^2$.  We apply this to the 
function $\phi=q'$ to see that 
$$
|q''(z)| \leq \frac{4}{(1-|z|^2)^2}\,.
$$
Finally, a result of Pommerenke ([20], p. 133) asserts that 
$$
(1-|z|^2)\left|\frac{h''(z)}{h'(z)}\right| \leq 2 + 
2\bigl(1+\tfrac12 \|{\Cal S}h\|\bigr)^{1/2}\,. \tag15
$$
Putting the estimates together, we conclude that 
$$
\|{\Cal S}f\| \leq \|{\Cal S}h\| +  
2\bigl(1+\tfrac12 \|{\Cal S}h\|\bigr)^{1/2} + 7\,.
$$

     Conversely, suppose that $\|{\Cal S}f\|<\infty$.  The formula (8) 
and the preceding estimates show that 
$$
|{\Cal S}h(z)| \leq |{\Cal S}f(z)| +  \frac{5}{(1-|z|^2)^2} + 
\frac{1}{1-|z|^2}\,\left|\frac{h''(z)}{h'(z)}\right|\,. \tag16
$$
In order to use Pommerenke's estimate of $h''/h'$, we apply the 
inequality (16) to the dilated function $f_r=h_r+\overline{g_r}$, 
where $0<r<1$ and $f_r(z)=f(rz)$.  Note that 
${\Cal S}f_r(z)=r^2{\Cal S}f(rz)$, so that 
$$
(1-|z|^2)^2|{\Cal S}f_r(z)| \leq (1-|rz|^2)^2|{\Cal S}f(rz)| 
\leq \|{\Cal S}f\|\,.
$$
Because $\|{\Cal S}h_r\|$ is finite for each $r<1$, we can apply (16) 
to $f_r$ and invoke (15) to infer that 
$$
(1-|z|^2)^2|{\Cal S}h_r(z)| \leq \|{\Cal S}f\| + 7 + 
2\bigl(1+\tfrac12 \|{\Cal S}h_r\|\bigr)^{1/2}\,,
$$
or
$$
\|{\Cal S}h_r\| - 2\bigl(1+\tfrac12 \|{\Cal S}h_r\|\bigr)^{1/2} \leq 
\|{\Cal S}f\| + 7\,.
$$
Now let $r\to1$ to conclude that $\|{\Cal S}h\|<\infty$\,.
\enddemo

     Finally, we observe that univalent harmonic mappings with range 
convex in one direction have finite Schwarzian norm.  These mappings are 
obtained by a known process of shearing conformal mappings whose range is 
convex in one direction.  (See [9], Section 3.4 for the shear 
construction.)  
\proclaim{Theorem 4} Suppose a function $\phi$ is analytic and univalent 
in the unit disk, and its range is convex in the horizontal direction.  
Let $f=h+\overline{g}$ be the harmonic shear of $\phi$ in the horizontal 
direction with dilatation $q^2$, where $q$ is analytic and $|q(z)|<1$ 
in the disk.  Then $\|{\Cal S}f\|<\infty$.  
\endproclaim
\demo{Proof} As shown in [3], the Schwarzian of $f$ has the form
$$
\align
{\Cal S}f=\ &{\Cal S}\phi+\frac{2({q'}^2+(1-q^2)qq'')}{(1-q^2)^2} 
-\frac{2qq'}{1-q^2}\frac{\phi''}{\phi'}\\
&+\frac{2\,\overline{q}}{1+|q|^2}\biggl\{q'' - 
q'\biggl(\frac{\phi''}{\phi'}+\frac{2qq'}{1-q^2}
\biggr)\biggr\} - 4\biggl(\frac{\overline{q}\,q'}
{1+|q|^2}\biggr)^2\, .
\endalign
$$
The preceding estimates for $q$ and its derivatives can now be 
applied to derive the inequality
$$
|{\Cal S}f(z)| \leq |{\Cal S}\phi(z)| + 
2\,\left|\frac{q''(z)}{1-q(z)^2}\right| + \frac{3}{1-|z|^2}\,
\left|\frac{\phi''(z)}{\phi'(z)}\right| + \frac{9}{(1-|z|^2)^2}\,.
$$
Since $|q'(z)/(1-q(z)^2)|\leq 1/(1-|z|^2)$, Cauchy's integral formula 
shows that the derivative satisfies 
$$
\left|\frac{q''(z)}{1 - q(z)^2} + \frac{2q(z)q'(z)^2}{(1 - q(z)^2)^2}
\right| \leq \frac{4}{(1-|z|^2)^2}\,,
$$
which implies that 
$$
\left|\frac{q''(z)}{1 - q(z)^2}\right| \leq \frac{6}{(1-|z|^2)^2}\,.
$$
Since $\phi$ is an analytic univalent function, a standard inequality 
({\it cf.} [8], p. 32) shows that  
$$
\left|\frac{\phi''(z)}{\phi'(z)}\right| \leq \frac{6}{1-|z|^2}\,,
$$
and Kraus' theorem gives $(1-|z|^2)^2|{\Cal S}\phi(z)| \leq 6$.  
Combining these estimates, we find that 
$$
\|{\Cal S}f\| \leq 6 + 12 + 18 + 9 = 45. 
$$
\enddemo

     The harmonic shear of the Koebe function with dilatation $z^2$ 
is $f=h+\overline{g}$, where 
$$
h(z) = \frac{\tfrac13}{(1-z)^3}\qquad \text{and} \qquad 
g(z) = \frac{z^2 - z + \tfrac13}{(1-z)^3} \,,
$$
up to additive constants.  Its Schwarzian is found to be 
$$
{\Cal S}f(z)=-4\biggl(\frac{1}{1-z}+\frac{\overline{z}}
{1+|z|^2}\biggr)^2.
$$
(See [3], p. 337, where the formula is printed incorrectly.)  
Note that $\|{\Cal S}f\|<\infty$, in agreement with Theorem 4.  
Although $h$ is not univalent, it has Schwarzian 
${\Cal S}h(z)=-4/(1-z)^2$ and Schwarzian norm $\|{\Cal S}h\|=16$.  
Since $\|{\Cal S}h\|<\infty$, the finiteness of $\|{\Cal S}f\|$ 
is also in agreement with Theorem 3. 

     {\it Addendum}.  In subsequent work [5] we have been able to show that 
the Schwarzian norms of univalent harmonic mappings are finite and have 
the uniform bound $\|{\Cal S}f\|<19,407$.

\enddocument